\documentclass[10pt,letterpaper,twocolumn]{article}

\usepackage[margin=0.75in]{geometry}
\usepackage[utf8]{inputenc}
\usepackage[T1]{fontenc}
\usepackage{lmodern}
\usepackage[numbers,sort&compress]{natbib}
\usepackage{tikz}
\usepackage{booktabs}
\usepackage{tabularx}
\usepackage{graphicx}
\usepackage{amsmath}
\usepackage{amssymb}
\usepackage{xcolor}
\usepackage{xurl}
\usepackage{caption}
\usepackage{microtype}
\usepackage{balance}
\usepackage{titlesec}
\usepackage[colorlinks=true,linkcolor=blue,citecolor=blue,urlcolor=blue]{hyperref}

\graphicspath{{media/}}
\usetikzlibrary{arrows.meta,positioning,fit,backgrounds}

\definecolor{PSEBlue}{HTML}{00509E}
\titleformat{\section}{\large\bfseries\color{PSEBlue}\raggedright}{}{0pt}{\MakeUppercase}
\titleformat{\subsection}{\normalsize\bfseries\color{PSEBlue}\raggedright}{}{0pt}{}
\titleformat{\subsubsection}{\normalsize\itshape\color{PSEBlue}}{}{0pt}{}
\titlespacing*{\section}{0pt}{12pt}{3pt}
\titlespacing*{\subsection}{0pt}{8pt}{3pt}

\newcommand{\PaperTitle}{}
\newcommand{\PaperAuthors}{}
\newcommand{\PaperAffiliations}{}
\newcommand{\PaperAbstract}{}
\newcommand{\PaperKeywords}{}

\definecolor{ORhfm}{HTML}{FDEBD0}
\definecolor{ORapm}{HTML}{FDEBD0}
\definecolor{ORmid}{HTML}{E67E22}
\definecolor{ORbot}{HTML}{FAC898}
\definecolor{ORbdr}{HTML}{A04000}
\definecolor{GRhfm}{HTML}{D5F5E3}
\definecolor{GRapm}{HTML}{D5F5E3}
\definecolor{GRmid}{HTML}{27AE60}
\definecolor{GRbot}{HTML}{A8DBA8}
\definecolor{GRbdr}{HTML}{1E8449}
\definecolor{GRtxt}{HTML}{145A32}

\renewcommand{\PaperTitle}{%
Towards Fewer Control Laws via Continuous-Time Multiparametric Programming}
\renewcommand{\PaperAuthors}{%
Lida Lamakani\textsuperscript{a,b} and
Efstratios N.~Pistikopoulos\textsuperscript{a,b*}}
\renewcommand{\PaperAffiliations}{%
\textsuperscript{a}Texas A\&M Energy Institute, Texas A\&M University,
College Station, TX, USA\\
\textsuperscript{b}Artie McFerrin Department of Chemical Engineering,
Texas A\&M University, College Station, TX, USA\\
\textsuperscript{*}Corresponding Author:
\href{mailto:stratos@tamu.edu}{stratos@tamu.edu}}
\renewcommand{\PaperAbstract}{%
Multiparametric programming offers a powerful solution to the computational burden of solving optimal control problems repeatedly online. By solving the problem once offline, it yields the optimal control laws as explicit closed-form functions of the initial system state, reducing online execution to a direct evaluation with no iterations required. Most existing works build this idea on a discrete-time foundation, slicing the time horizon into intervals and applying KKT conditions to the resulting algebraic system. This discretization forces a tradeoff: too few intervals and the model fails to capture the true system dynamics, while too many cause the problem size, the number of decision variables, and the number of critical regions to grow rapidly, making both offline preparation and online lookup increasingly expensive.
This work develops a multiparametric framework that works directly with the continuous-time problem. Pontryagin's Maximum Principle (PMP) is applied without any model discretization, and the optimal control is recovered as an explicit function of the initial state. Compared to the discrete-time formulation, the continuous-time approach produces substantially fewer critical regions, and this number remains fixed regardless of accuracy requirements, since it reflects the structure of the problem itself rather than a discretization grid. The framework also yields the switching times as explicit functions of the initial state, directly exposing when and how the optimal control structure changes over the horizon. Knowing these switching times in advance allows the real-time controller to skip unnecessary computations between them, further reducing the online execution cost. Results from a PAROC framework case study demonstrate that the continuous-time multiparametric approach is a rigorous alternative to the conventional discrete-time formulations.}
\renewcommand{\PaperKeywords}{%
Multiparametric dynamic optimization, Discretization-free control, Continuous-time switching instants}

\newcommand{\arxivmaketitle}{%
  \twocolumn[{%
    \begin{center}
      {\LARGE\bfseries \PaperTitle\par}
      \vspace{0.75em}
      {\large \PaperAuthors\par}
      \vspace{0.5em}
      {\small \PaperAffiliations\par}
      \vspace{0.5em}
      \vspace{0.75em}
      \begin{minipage}{0.92\textwidth}
        \begin{abstract}
        \PaperAbstract
        \end{abstract}
        \noindent\textbf{Keywords:} \PaperKeywords
      \end{minipage}
    \end{center}
    \vspace{1em}
  }]%
}

\begin{document}

\arxivmaketitle

\section{Introduction}

In multiparametric model predictive control, the optimal control
problem is solved once offline as an explicit function of the
initial state $x_0$, dividing the parameter space into critical
regions, each carrying its own control law. Online execution
reduces to identifying the critical region for the current
state and evaluating the corresponding law. The number of critical
regions sets the memory required by the controller and the time
required to identify the active region at each sampling
instant.~\cite{pistikopoulos2020multi}

The PAROC framework~\cite{pistikopoulos2015paroc}
applies multiparametric MPC to chemical process systems. It
reduces a high-fidelity plant to a linear surrogate, solves the
optimal control problem offline as a function of the initial
state, and deploys the resulting set of critical regions online,
with the controller identifying its current critical region and
evaluating the corresponding control law without any online
optimization. The 32-state Newell-Lee distillation
column~\cite{newell1989applied} has served as a standard benchmark
for this framework~\cite{pistikopoulos2015paroc} and is used here as
the benchmark setting for comparing the continuous-time formulation
with two discrete-time multiparametric counterparts.

Physical systems evolve continuously in time. When a discrete-time
(DT) surrogate is required, two identification routes are
available: identify the ODE
$\dot{x}=A_\mathrm{c}x+B_\mathrm{c}u$ first and then discretize it
(DT-ODE)~\cite{pistikopoulos2007multi}, or fit the algebraic model
$x_{k+1}=A_d x_k+B_d u_k$ directly from sampled input-output
data (DT-direct)~\cite{pistikopoulos2015paroc}. The classical PAROC
framework uses the latter~\cite{pistikopoulos2015paroc}. In both
routes, the active-set structure depends on the chosen time grid.
Finer grids introduce more decision stages and, therefore, more
possible active-set changes in the explicit solution.

This paper develops a continuous-time (CT) multiparametric
control-law construction for reduced-order process models and
demonstrates it within the PAROC workflow on the 32-state
Newell-Lee distillation column~\cite{pistikopoulos2015paroc}.

A CT ODE surrogate and two DT surrogates (DT-direct and DT-ODE)
are identified from the same 32-state simulation data and validated
against the original nonlinear column. The CT multiparametric
procedure~\cite{lamakani2026multiparametric} returns a
critical-region map for this case study by expressing the state and
costate trajectories through Hamiltonian matrix exponentials. In the
CT map, switching times divide the horizon into time segments with
fixed active constraint sets, and these switching times depend on
the initial state. In the DT maps, active constraints are assigned
at the $N$ sampling times, creating more possible active-set
combinations and, as a result, both DT formulations return more
critical regions.

\section{Classical and Next-Generation PAROC Frameworks}

Figs.~\ref{fig:paroc_classic} and~\ref{fig:paroc_ct} show the
two frameworks side by side. The high-fidelity modeling stage and
the online controller are identical in both. The difference lies
entirely in the offline preparation: what model is identified in
Stage~2 and what optimality conditions are applied in Stage~3.

\begin{figure}[!t]
\centering
\begin{minipage}[b]{0.85\columnwidth}
\centering
\resizebox{\linewidth}{!}{%
\begin{tikzpicture}[
  node distance = 0.7cm,
  every node/.style = {align=center,
    font=\fontsize{10}{12}\selectfont, text=black},
  HFM/.style = {draw=ORbdr, fill=ORhfm, rounded corners=3pt,
                line width=0.9pt, text width=8cm,
                minimum height=1.2cm, inner sep=7pt},
  SID/.style = {HFM}, APR/.style = {HFM},
  BOT/.style = {draw=ORbdr, fill=ORbot, rounded corners=3pt,
                line width=0.9pt, text width=8cm,
                minimum height=1.2cm, inner sep=7pt},
  MID/.style = {draw=ORbdr, fill=ORmid, rounded corners=3pt,
                line width=0.9pt, text width=8cm,
                minimum height=1.2cm, inner sep=7pt},
  ARW/.style = {-{Stealth[length=8pt,width=5pt]},
                ORbdr, line width=1.1pt}
]
\node[HFM] (hfm)
  {Process `High-Fidelity' Dynamic Modeling};
\node[SID, below=of hfm] (sid)
  {System Identification \& Model Reduction};
\node[APR, below=of sid] (apr)
  {Approximate Model --- \textbf{in Algebraic Space}\\[5pt]
   $x_{k+1} = A_d\,x_k + B_d\,u_k$};
\node[BOT, below=of apr] (mpp)
  {DT Multiparametric Programming\\[5pt]
   KKT optimality conditions};
\node[MID, below=of mpp] (ctl)
  {Multi-parametric Control \& Receding Horizon Policies};
\draw[ARW](hfm)--(sid);
\draw[ARW](sid)--(apr);
\draw[ARW](apr)--(mpp);
\draw[ARW](mpp)--(ctl);
\begin{scope}[on background layer]
  \node[draw=ORbdr, fill=orange!5, dashed, rounded corners=6pt,
        line width=0.8pt, fit=(sid)(apr)(mpp),
        inner sep=13pt] (ob){};
\end{scope}
\node[font=\small\bfseries\itshape, text=ORbdr,
      anchor=south east, inner sep=2.2pt]
      at (ob.south east) {Offline, once};
\end{tikzpicture}
}%
\caption{Classical PAROC framework. Number of critical regions grows with $N$.}
\label{fig:paroc_classic}
\end{minipage}\hfill
\begin{minipage}[b]{0.85\columnwidth}
\centering
\resizebox{\linewidth}{!}{%
\begin{tikzpicture}[
  node distance = 0.7cm,
  every node/.style = {align=center,
    font=\fontsize{10}{12}\selectfont, text=black},
  HFM/.style = {draw=GRbdr, fill=GRhfm, rounded corners=3pt,
                line width=0.9pt, text width=8cm,
                minimum height=1.2cm, inner sep=7pt},
  SID/.style = {HFM}, APR/.style = {HFM},
  BOT/.style = {draw=GRbdr, fill=GRbot, rounded corners=3pt,
                line width=0.9pt, text width=8cm,
                minimum height=1.2cm, inner sep=7pt},
  MID/.style = {draw=GRbdr, fill=GRmid, rounded corners=3pt,
                line width=0.9pt, text width=8cm,
                minimum height=1.2cm, inner sep=7pt},
  ARW/.style = {-{Stealth[length=8pt,width=5pt]},
                GRbdr, line width=1.1pt}
]
\node[HFM] (hfm)
  {Process `High-Fidelity' Dynamic Modeling};
\node[SID, below=of hfm] (sid)
  {System Identification \& Model Reduction};
\node[APR, below=of sid] (apr)
  {Approximate Model --- \textbf{in ODE Space}\\[5pt]
   $\dot{x} = A_c\,x + B_c\,u$};
\node[BOT, below=of apr] (mpp)
  {CT Multiparametric Programming\\[5pt]
   PMP optimality conditions};
\node[MID, below=of mpp] (ctl)
  {Multi-parametric Control \& Receding Horizon Policies};
\draw[ARW](hfm)--(sid);
\draw[ARW](sid)--(apr);
\draw[ARW](apr)--(mpp);
\draw[ARW](mpp)--(ctl);
\begin{scope}[on background layer]
  \node[draw=GRbdr, fill=green!5, dashed, rounded corners=6pt,
        line width=0.8pt, fit=(sid)(apr)(mpp),
        inner sep=14pt] (ob){};
\end{scope}
\node[font=\small\bfseries\itshape, text=GRbdr,
      anchor=south east, inner sep=2.2pt]
      at (ob.south east) {Offline, once};
\end{tikzpicture}
}%
\caption{Next-generation PAROC framework (this work). Critical regions fixed by switching structure.}
\label{fig:paroc_ct}
\end{minipage}
\end{figure}

In the classical framework (Fig.~\ref{fig:paroc_classic}), Stage~2
produces $x_{k+1}=A_d x_k+B_d u_k$ and Stage~3 applies KKT
conditions to this algebraic model, returning the piecewise-affine
explicit solution~\cite{bemporad2002explicit,tondel2003algorithm}.

In the next-generation framework (Fig.~\ref{fig:paroc_ct}),
Stage~2 fits the vector field $\dot{x}=A_\mathrm{c}x+B_\mathrm{c}u$
from state-derivative triples $(\dot{x},x,u)$. This ODE enters
Stage~3 directly, without discretization. PMP characterizes the
optimal input through arc structures, and the critical region map
reflects those structures rather than the time
grid~\cite{sakizlis2005explicit}.

\section{Continuous-Time (CT) Multiparametric Optimal Control}

The full theory, PMP conditions, and proof that critical region
boundaries are hyperplanes for LQ systems with input bounds are
in~\cite{lamakani2026linear}. This section states the optimal
control problem formulation and the key PMP results applied in
the case study that follows.
\begin{equation}
  \begin{aligned}
    \min_{u(\cdot)}\quad & \tfrac{1}{2}x(t_f)^\top P_{\!f}x(t_f)
    +\int_{t_0}^{t_f}\tfrac{1}{2}\bigl(x(t)^\top Qx(t)+Ru(t)^2\bigr)\,dt\\
    \mathrm{s.t.}\quad & \dot{x}(t)=Ax(t)+Bu(t), \quad t\in[t_0,t_f],\\
    & x(t_0)=\theta,\\
    & |u(t)|\leq u_\mathrm{max}.
  \end{aligned}
  \label{eq:ocp}
\end{equation}
PMP converts this infinite-dimensional optimization into a boundary
value problem by introducing the costate $\lambda(t)$ and forming
the Hamiltonian $H=L+\lambda^\top(Ax+Bu)+\mu_U(u-u_\mathrm{max})
+\mu_L(-u-u_\mathrm{max})$. The necessary conditions
$\dot{x}=\nabla_\lambda H$, $\dot{\lambda}=-\nabla_x H$,
$\nabla_u H=0$, $\lambda(t_f)=P_{\!f}x(t_f)$, together with
corner conditions at switching times, fully characterize the
optimal solution~\cite{lamakani2026multiparametric}.

The optimal solution over $[t_0,t_f]$ consists of a sequence of
arcs, where switching times $t_s\in(t_0,t_f)$ mark the
transitions between consecutive arcs. On each arc, the set of
active constraints remains fixed, and the ODEs governing
$x$, $\lambda$, $\mu$, and $g$ are linear. Their solutions,
together with the optimal input $u$, are exponential functions
of the arc's initial state, $t_s$, and $t$, computed analytically
from matrix exponentials of the Hamiltonian. Since the initial
state of each arc is determined by $\theta$ through the preceding
arcs, $x$, $u$, $\lambda$, $\mu$, and $g$ are ultimately
exponential functions of $\theta$, $t_s$, and $t$. For brevity,
these arguments are suppressed in the remainder of this section.
The active set on each arc is defined as
$\mathcal{A}=\{i:g_i=0\}$, and a critical region is the set of
all initial conditions $\theta\in\Theta$ whose optimal solution
follows the same arc sequence, that is, the same ordered sequence
of active sets $\mathcal{A}$ from $t_0$ to $t_f$.

The critical region boundaries are the values of $\theta$ where
the arc sequence changes. A constraint becomes active at the
value of $\theta$ where $\min_t \mu=0$ on the corresponding
constrained arc, and becomes inactive where $\max_t g=0$ on the
corresponding unconstrained arc. Since the costate is linear in
$\theta$ at any fixed time, and the boundary-defining extrema of
$\mu$ and $g$ are attained at the fixed horizon endpoints $t_0$
and $t_f$ for this system~\cite{lamakani2026linear}, every
critical region boundary is a hyperplane in $\Theta$. These
hyperplanes divide $\Theta$ into polyhedral critical regions,
each carrying a unique arc sequence. Within each region, the
optimal control law is recovered directly from the PMP conditions
as a function of $\theta$, and online execution requires nothing
beyond locating the current state in its region and evaluating
that law.

\section{Binary Distillation Column}

\subsection{High-Fidelity Model}

The 32-tray Newell-Lee column~\cite{newell1989applied} governs tray
compositions $x_i$ through:
\begin{align}
  \dot{x}_1    &= \tfrac{1}{A_\mathrm{cond}}V(y_2-x_1),
                  \label{eq:cond}\\
  \dot{x}_i    &= \tfrac{1}{A_\mathrm{tray}}
                  \bigl[L_1(x_{i-1}-x_i)-V(y_i-y_{i+1})\bigr],
                  i=2,\ldots,16,\label{eq:rect}\\
  \dot{x}_{17} &= \tfrac{1}{A_\mathrm{tray}}
                  \bigl[Fx_F+L_1x_{16}-L_2x_{17}
                  -V(y_{17}-y_{18})\bigr],\label{eq:feed}\\
  \dot{x}_i    &= \tfrac{1}{A_\mathrm{tray}}
                  \bigl[L_2(x_{i-1}-x_i)-V(y_i-y_{i+1})\bigr],
                i=18,\ldots,31,\label{eq:strip}\\
  \dot{x}_{32} &= \tfrac{1}{A_\mathrm{reb}}
                  \bigl[L_2x_{31}-(F-D)x_{32}-Vy_{32}\bigr],
                  \label{eq:reb}
\end{align}
with $y_i=\alpha x_i/(1+(\alpha-1)x_i)$, $V=L_1+D$, $L_2=F+L_1$.
Eq.~\eqref{eq:cond} is the total condenser,
Eqs.~\eqref{eq:rect} the rectifying section,
Eq.~\eqref{eq:feed} the feed tray, Eqs.~\eqref{eq:strip}
the stripping section, and Eq.~\eqref{eq:reb} the reboiler.
Column parameters~\cite{pistikopoulos2015paroc} are $F=0.40$,
$D=0.20$, $A_\mathrm{cond}=0.50$, $A_\mathrm{tray}=0.25$,
$A_\mathrm{reb}=1.00$, $\alpha=1.60$, with the reflux-ratio
deviation $u=\mathrm{RR}-2.70$ as the single manipulated variable.
The reduced state tracks deviations at the two product trays:
$x=[x_1-x_1^\mathrm{nom},\,x_{32}-x_{32}^\mathrm{nom}]^\top$,
where $x_1^\mathrm{nom}$ and $x_{32}^\mathrm{nom}$ are the
steady-state top and bottom compositions of the 32-state column
at the nominal operating point $\mathrm{RR}=2.70$, $x_F=0.50$.
The initial condition $\theta=x(t_0)\in\Theta\subset\mathbb{R}^2$
is the multiparametric parameter in all three formulations; its
two components are the top and bottom composition deviations at
the start of each control horizon.

All computations were carried out in Python: SciPy was used for
nonlinear simulation, least-squares identification, and the CT
PMP-based computations, including Hamiltonian matrix exponentials
and switching-time root finding, while PPOPT with the combinatorial
algorithm~\cite{kenefake2022ppopt,pistikopoulos2020multi} was used
for the two DT mp-QPs.

\subsection{Surrogate Identification}

All identifications use 36 closed-loop simulation trajectories
of the 32-state column under piecewise-constant reflux excitations.

For the CT surrogate, regression on state-derivative triples
$(\dot{x}(t),x(t),u(t))$ yields $\dot{x}=A_c x+B_cu$:
\begin{equation}
  A_c=
  \begin{bmatrix}-3.795&-3.941\\1.821&1.873\end{bmatrix},\quad
  B_c=
  \begin{bmatrix}4.35\!\times\!10^{-4}\\
  -2.17\!\times\!10^{-3}\end{bmatrix},
  \label{eq:AB_ct}
\end{equation}
with derivative-fit RMSE $=2.93\!\times\!10^{-5}$ and $R^2=0.993$.

For the DT-direct surrogate, regression on sampled state
pairs $(x_{k+1},x_k,u_k)$ at $h=0.1\,\mathrm{min}=6\,\mathrm{s}$
yields
\[
x_{k+1}=A_d x_k+B_d u_k:
\]
\begin{equation}
A_d=\begin{bmatrix}
0.647 & -0.366\\
0.163 &  1.168
\end{bmatrix},\quad
B_d=\begin{bmatrix}
5.54\times 10^{-5}\\
-2.24\times 10^{-4}
\end{bmatrix},
\label{eq:AB_dt}
\end{equation}
with one-step RMSE $=1.12\!\times\!10^{-6}$ and $R^2=0.999$.

For the DT-ODE surrogate, the CT ODE surrogate is
discretized by exact zero-order hold at
$h=0.1\,\mathrm{min}=6\,\mathrm{s}$:
\begin{equation}
  A_d^\mathrm{ODE}=
  \begin{bmatrix}0.654 & -0.358\\0.166 &  1.170\end{bmatrix},\quad
  B_d^\mathrm{ODE}=
  \begin{bmatrix}7.58\!\times\!10^{-5}\\
  -2.32\!\times\!10^{-4}\end{bmatrix},
  \label{eq:AB_zoh}
\end{equation}
with the stage cost discretized as
$\tfrac{1}{2}\sum_{k=0}^{N-1}h(x_k^\top Qx_k+Ru_k^2)
+\tfrac{1}{2}x_N^\top P_{\!f}x_N$.

\begin{figure}[t]
  \centering
  \includegraphics[width=1\columnwidth]{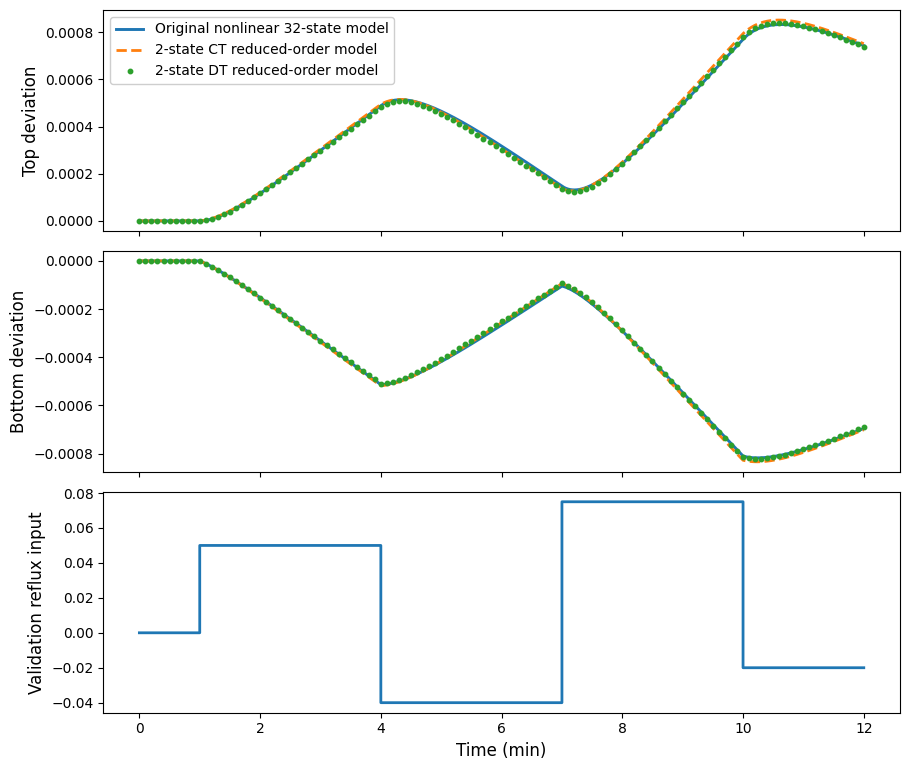}
  \caption{Validation of the two reduced-order surrogates against the original nonlinear 32-state column.}
  \label{fig:val}
\end{figure}

The use of a two-state surrogate is supported by the Hankel
singular values (HSVs) of the 32-state input-output linearization.
Table~\ref{tab:hsv} shows that two dominant dynamic directions
capture 98.50\% of the cumulative input-output energy, with a
balanced-truncation error bound of $1.94\times 10^{-3}$ for the
discarded dynamics. These HSVs support the low-order nature of the
input-output behavior, while the selected reduced states are the
top and bottom composition deviations defined above.
Fig.~\ref{fig:val} validates the reduced models against the
original nonlinear 32-state column using a separate stepwise input
trajectory covering the allowed reflux range $u\in[-0.08,0.08]$.
At the DT sampling times, the combined validation RMSEs are
$1.41\times 10^{-5}$ for the CT surrogate and
$1.47\times 10^{-5}$ for the DT-direct surrogate.

\begin{table}[t]
  \caption{HSVs of the 32-state input-output linearization; two dominant dynamic directions capture 98.50 percent of the cumulative input-output energy.}
  \label{tab:hsv}
  \centering
  \renewcommand{\arraystretch}{1.10}
  \begin{tabularx}{\columnwidth}{@{}cXX@{}}
    \toprule
    \textbf{Mode} & \textbf{HSV} &
    \textbf{Cumulative fraction}\\
    \midrule
    1 & $6.23\!\times\!10^{-2}$ & 96.10\%\\
    2 & $1.55\!\times\!10^{-3}$ & 98.50\%\\
    3 & $7.20\!\times\!10^{-4}$ & 99.61\%\\
    $\vdots$ & $\vdots$ & $\vdots$\\
    32 & -- & 100\%\\
    \bottomrule
  \end{tabularx}
\end{table}

\subsection{CT and DT Multiparametric Solutions}
\label{sec:comp}

All three formulations use identical objectives: $Q=\mathrm{diag}(1500,300)$, $R=1$, $P_{\!f}=2.5Q$, horizon $t_{\!f}=1$ min,
$u_\mathrm{max}=0.08$, $\Theta=[-0.020,0.020]\times[-0.010,0.010]$.

\begin{figure}[t]
  \centering
  \includegraphics[width=1\columnwidth]{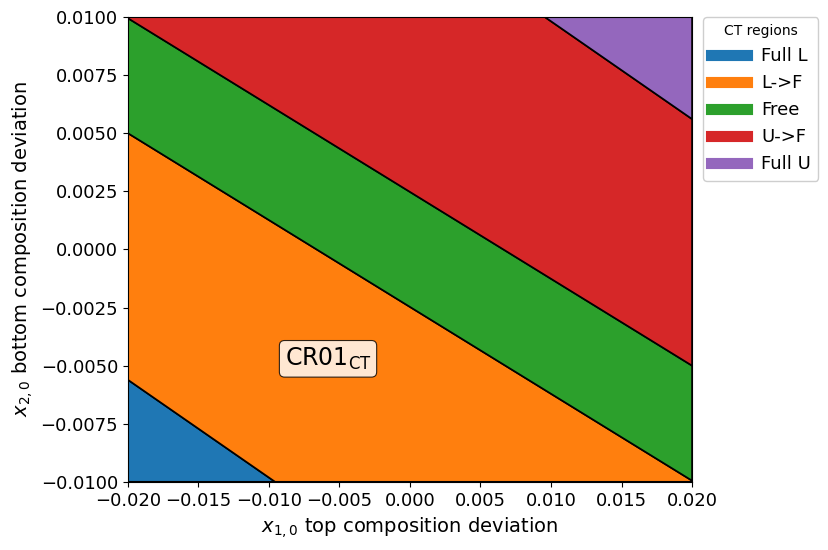}
  \caption{CT solution: 5 critical regions. Boundaries from Hamiltonian matrix exponentials.}
  \label{fig:ct}
\end{figure}

Following the procedures of~\cite{lamakani2026linear,lamakani2026multiparametric},
the CT multiparametric solution for this problem yields five critical
regions with hyperplane boundaries computed from the Hamiltonian
matrix exponential (Fig.~\ref{fig:ct}). Each region corresponds to
a distinct arc sequence: in the Free region the input is
unconstrained over the entire horizon; in the single-switch
regions Upper-to-Free and Lower-to-Free the input is initially
saturated at the upper or lower bound and switches to the
unconstrained arc at some $t_s\in(t_0,t_f)$; and in the fully
constrained regions Full Upper and Full Lower the input remains
saturated at the upper or lower bound throughout.

\begin{figure}[t]
\centering
\begin{minipage}[b]{1\columnwidth}
\centering
\includegraphics[width=\linewidth]{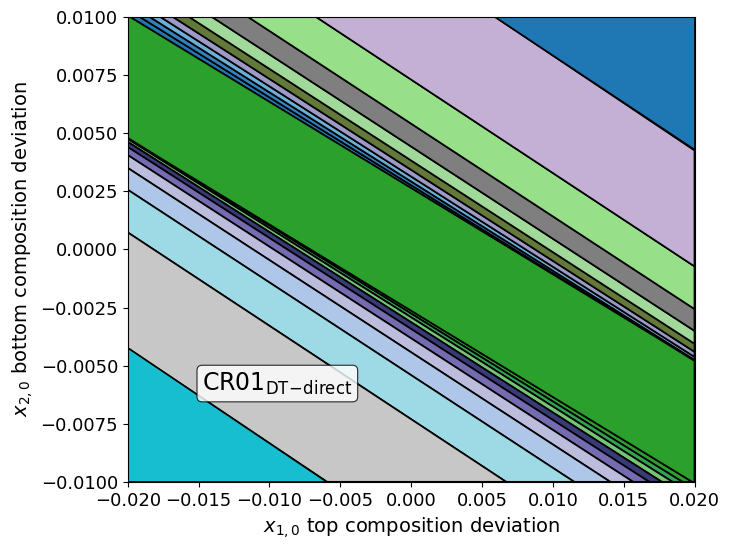}
\caption{DT-direct route ($N=10$): 23 CRs.}
\label{fig:dt_direct}
\end{minipage}\hfill
\begin{minipage}[b]{1\columnwidth}
\centering
\includegraphics[width=\linewidth]{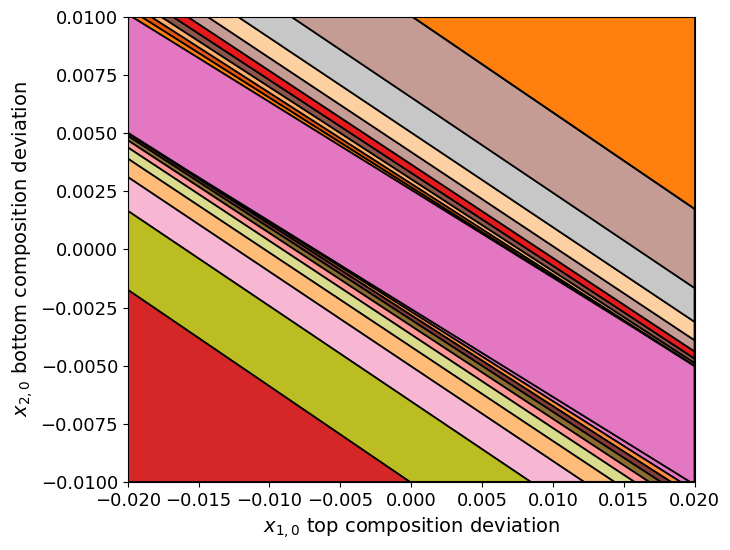}
\caption{DT-ODE route ($N=10$): 23 CRs.}
\label{fig:dt_zoh}
\end{minipage}
\end{figure}

Both DT formulations are solved by the multiparametric
combinatorial algorithm at $N=10$, returning 23 critical regions
each (Figs.~\ref{fig:dt_direct},~\ref{fig:dt_zoh}). The larger
number of DT regions compared with the CT solution reflects the
possible active-set combinations over the $N$ discrete time steps.
Although both routes yield the same number of critical regions in
this case, their boundaries differ because the DT-direct and DT-ODE
surrogates approximate the same plant differently, and neither is
systematically more accurate than the other.

Table~\ref{tab:cr01_laws} reports the explicit control laws for the two representative regions labeled in the maps: $\mathrm{CR01}_{\mathrm{CT}}$ in Fig.~\ref{fig:ct} and $\mathrm{CR01}_{\mathrm{DT\text{-}direct}}$ in Fig.~\ref{fig:dt_direct}. The table lists their arc/active-set descriptions, critical-region boundaries, and input laws.

\begin{table}[htp!]
\caption{Explicit control laws in $\mathrm{CR01_{CT}}$ and $\mathrm{CR01_{DT-direct}}$.}
\label{tab:cr01_laws}
\centering\scriptsize
\renewcommand{\arraystretch}{1.20}
\setlength{\tabcolsep}{1pt}
\begin{tabularx}{\columnwidth}{@{}l>{\raggedright\arraybackslash}p{0.50\columnwidth}>{\raggedright\arraybackslash}p{0.33\columnwidth}@{}}
\toprule
 & $\mathrm{CR01_{CT}}$ & $\mathrm{CR01_{DT-direct}}$ \\
\midrule
Arc/set
  & $L\to F$
  & $[\underbrace{-1,\cdots,-1}_{9},0]$ \\[3pt]
Boundary
  & $\begin{bmatrix}-0.3878&-0.9217\\0.3502&0.9367\end{bmatrix}\theta
    \le\begin{bmatrix}0.01292\\-0.002317\end{bmatrix}$
  & $\begin{bmatrix}0.3728&0.9279\\-0.3781&-0.9258\end{bmatrix}\theta
    \le\begin{bmatrix}-6.78{\times}10^{-3}\\1.15{\times}10^{-2}\end{bmatrix}$
    \\[3pt]
Input
  & $u(t)=-u_\mathrm{max},\quad t_0\!\le\!t\!\le\!t_s$\newline
    $u(t)=-R^{-1}B_c^\top\lambda(t),\quad t_s\!<\!t\!\le\!t_f$\newline
    $[x(t);\lambda(t)]^\top\!=\!e^{H_F(t-t_s)}[x_s;\,P_{\!f}x_s]^\top$\newline
    $x_s\!=\!e^{A_c t_s}\theta\!+\!(e^{A_c t_s}\!-\!I)A_c^{-1}B_c(-u_\mathrm{max})$
  & $u_k=-u_\mathrm{max},\quad k\!=\!0,\ldots,8$\newline
    $u_9\!=\!2.654\,\theta_1\!+\!6.500\,\theta_2$\newline
    $\phantom{u_9}+\,7.560{\times}10^{-4}$ \\[3pt]
$t_s(\theta)$ & see~\eqref{eq:ts_fit} & -- \\
\bottomrule
\end{tabularx}
\end{table}

$H_L$ and $H_F$ are the $4{\times}4$ PMP Hamiltonian matrices:
\begin{equation}
H_L=\begin{bmatrix}A_c & 0\\-Q & -A_c^\top\end{bmatrix},\qquad
H_F=\begin{bmatrix}A_c & -B_cR^{-1}B_c^\top\\-Q & -A_c^\top\end{bmatrix}.
\label{eq:hamiltonians}
\end{equation}
The switching time for $\mathrm{CR01_{CT}}$, with
$\bar{x}_i=\theta_i/\theta_{i,\max}$:
\begin{equation}
t_s=\dfrac{N_s}{D_s}\ \mathrm{[sec]},
\label{eq:ts_fit}
\end{equation}
\footnotesize
\begin{flalign*}
N_s &= 133.1+917.9\bar{x}_1+1208.6\bar{x}_2
       +1557.2\bar{x}_1^2+4112.3\bar{x}_1\bar{x}_2 &\\
    &\quad+2708.3\bar{x}_2^2-9.33\bar{x}_1^3
       -9.65\bar{x}_1^2\bar{x}_2-4.96\bar{x}_1\bar{x}_2^2
       -9.85\bar{x}_2^3, &\\
D_s &= 1+10.48\bar{x}_1+13.44\bar{x}_2
       +23.69\bar{x}_1^2+62.19\bar{x}_1\bar{x}_2
       +41.14\bar{x}_2^2, &
\end{flalign*}
\normalsize
with $R^2=0.99999$.

\section{Results and Discussion}

The CT formulation yields 5 critical regions for this case, each
corresponding to a distinct optimal arc sequence. This number is
fixed for a given horizon and reflects the switching structures
present in $\Theta$. Both the DT-direct and DT-ODE routes at $N=10$
yield 23 regions. In the DT solutions, the parameter space is divided
into more regions, each associated with a specific combination of
active constraints at the discrete time steps, which makes the
physical interpretation of the control strategy less direct.

Table~\ref{tab:u0} reports $u_0$ from all three formulations at
ten representative initial conditions. The percentages are relative
deviations from the CT first control move, computed as
$100(u_0^{\mathrm{DT}}-u_0^{\mathrm{CT}})/u_0^{\mathrm{CT}}$.
They quantify differences between the CT and DT control laws, not
surrogate-identification errors. For constrained regions, all three
agree exactly. For points in the Free region, both DT routes introduce
deviations relative to the CT solution, as shown in the table. Since
the switching time $t_s(\theta)$ varies continuously with the initial
condition, no fixed set of sampling instants can capture it exactly
for all $\theta\in\Theta$. Removing this error would require making
the sampling times optimization variables, which introduces
bilinear terms into the KKT conditions and eliminates the quadratic
programming structure required for a piecewise-affine explicit
solution. This approximation error is therefore an inherent limitation
of DT multiparametric programming, independent of the particular
discretization used.

\begin{table}[!htp]
\caption{First control move from CT, DT-direct, and DT-ODE at N=10. $\Delta_d(\%)$ and $\Delta_z(\%)$ are relative deviations from the CT first control move.}
  \label{tab:u0}
  \centering
  \renewcommand{\arraystretch}{1.05}
  \setlength{\tabcolsep}{3pt}
  \footnotesize
  \begin{tabular}{@{}rr rrr rr@{}}
    \toprule
    $\theta_1$ & $\theta_2$ &
    $u_0^{\mathrm{CT}}$ &
    $u_0^{\mathrm{DT-d}}$ & $\Delta_d(\%)$ &
    $u_0^{\mathrm{DT-z}}$ & $\Delta_z(\%)$\\
    \midrule
    $ 0.0050$ & $-0.0050$ & $-0.0800$ & $-0.0800$ & $ 0.00$ & $-0.0800$ & $ 0.00$\\
    $-0.0175$ & $ 0.0100$ & $ 0.0800$ & $ 0.0800$ & $ 0.00$ & $ 0.0800$ & $ 0.00$\\
    $ 0.0050$ & $ 0.0100$ & $ 0.0800$ & $ 0.0800$ & $ 0.00$ & $ 0.0800$ & $ 0.00$\\
    $ 0.0076$ & $-0.0042$ & $-0.0439$ & $-0.0419$ & $-4.73$ & $-0.0415$ & $-5.53$\\
    $-0.0020$ & $ 0.0017$ & $ 0.0308$ & $ 0.0291$ & $-5.65$ & $ 0.0295$ & $-4.21$\\
    $ 0.0048$ & $-0.0025$ & $-0.0228$ & $-0.0218$ & $-4.41$ & $-0.0214$ & $-5.99$\\
    $-0.0011$ & $ 0.0009$ & $ 0.0166$ & $ 0.0156$ & $-5.64$ & $ 0.0159$ & $-4.23$\\
    $-0.0010$ & $ 0.0017$ & $ 0.0429$ & $ 0.0403$ & $-6.00$ & $ 0.0413$ & $-3.70$\\
    $ 0.0010$ & $-0.0018$ & $-0.0461$ & $-0.0433$ & $-6.02$ & $-0.0444$ & $-3.68$\\
    $-0.0150$ & $-0.0050$ & $-0.0800$ & $-0.0800$ & $ 0.00$ & $-0.0800$ & $ 0.00$\\
    \bottomrule
  \end{tabular}
\end{table}

\begin{figure}[t]
\centering
\begin{minipage}[b]{1\columnwidth}
\centering
\includegraphics[width=\linewidth]{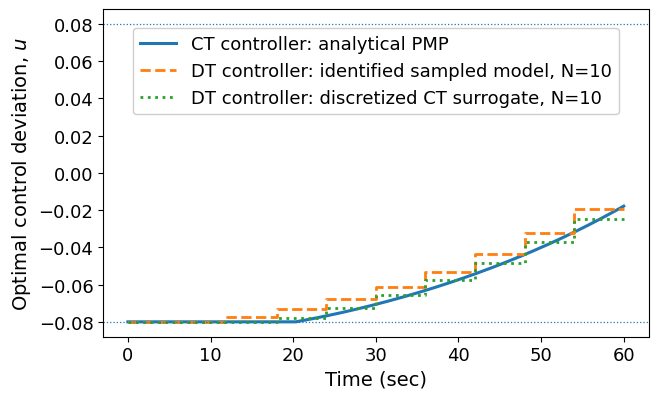}
\caption{theta=(-0.01, 0.001): CT switches at t=20.32 sec; DT-direct route at t=12.00 sec; DT-ODE route at t=18.00 sec.}
\label{fig:traj1}
\end{minipage}\hfill
\begin{minipage}[b]{1\columnwidth}
\centering
\includegraphics[width=\linewidth]{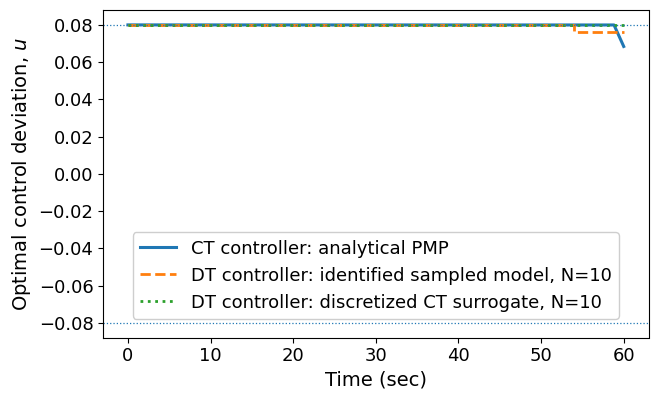}
\caption{theta=(0.015, 0.0057): CT switches at t=59.16 sec; DT-direct route at t=54.00 sec; DT-ODE route stays constrained.}
\label{fig:traj2}
\end{minipage}
\end{figure}

Figs.~\ref{fig:traj1} and~\ref{fig:traj2} compare optimal
trajectories at two initial conditions. For
$\theta=(-0.01,0.001)$, all three agree on $u_0=-0.08$; both DT
routes approximate the smooth CT trajectory with a piecewise-constant
input, switching at the nearest grid point rather than at the exact
CT switching time. For $\theta=(0.015,0.0057)$, the CT solution
switches from Full Upper to Free at $t=59.15515\,\mathrm{sec}$, the
DT-direct route switches at $t=54.00\,\mathrm{sec}$, and the DT-ODE
route remains constrained for the entire horizon. The CT formulation
provides the switching time $t_s(\theta)$ explicitly as a function of
the initial condition, specifying when the optimal control strategy
changes and accelerating online implementation by using the
precomputed switching time to know whether a sampling update occurs
before or after the switch.

\section{Conclusions}

In a CT multiparametric solution, critical regions correspond to
distinct arc sequences of the optimal input, and their number is
fixed for a given horizon. In a DT solution, critical regions
correspond to active-set combinations at each of the $N$ time steps,
and their number can grow with $N$. On the 32-state Newell-Lee
distillation column, the CT formulation returned 5 critical regions,
while both the DT-direct and DT-ODE routes returned 23 critical
regions at $N=10$. All three formulations yielded the same $u_0$ in
the constrained regions. In the Free region, both DT routes introduced
deviations relative to the CT solution, and the DT-ODE route failed
to identify the correct arc sequence at an initial condition near a
switching boundary, despite starting from the same CT ODE surrogate.

Integrating the CT framework into PAROC requires changing only the
surrogate identification and the multiparametric solver;
high-fidelity modeling, closed-loop validation, and online
deployment are unchanged. Future work will extend the framework to
problems where critical region boundaries become curved and require
continuation-based methods, and to robust formulations that handle
bounded plant uncertainty conditions.

\section{Acknowledgements}

\section{Declaration of Use of AI}

AI-assisted writing tools were used for editorial support during
manuscript preparation. All technical content, mathematical
derivations, numerical results, and conclusions are the sole work
of the authors.

\section{Author Identifiers}

Author ORCIDs:

Lamakani L: 0009-0006-5142-194X

Pistikopoulos EN: 0000-0001-6220-818X

\balance

\end{document}